\documentclass[11pt]{article}
\usepackage{amsmath}
\usepackage{amsfonts}
\usepackage{graphicx}
\usepackage{enumerate}
\usepackage{natbib}
\usepackage[papersize={8.5in,11in},margin=1in]{geometry}
\usepackage{lineno}

\begin{document}
\begin{center}
\Large\bf{What the collapse of the ensemble Kalman filter\\ tells us about particle filters}
\end{center}

\begin{center}
  Matthias Morzfeld  \\
  \emph{Department of Mathematics, University of Arizona}
  \vspace{3mm}
  
  Daniel Hodyss\\
  \emph{Naval Research Laboratory}
  \vspace{3mm}
  
  Chris Snyder\\
  \emph{National Center for Atmospheric Research}
  \vspace{6mm}
\end{center}

\begin{center}
\emph{Abstract}
\end{center}
\vspace{-3mm}
The ensemble Kalman filter (EnKF) is a reliable data assimilation tool 
for high-dimensional meteorological problems.
On the other hand, the EnKF can be interpreted as a particle filter,
and particle filters collapse in high-dimensional problems.
We explain that these seemingly contradictory statements
offer insights about how particle filters function in certain high-dimensional problems,
and in particular support recent efforts in meteorology to ``localize'' particle filters,
i.e., to restrict the influence of an observation to its neighborhood.

\section{Introduction}
Many applications in geophysics 
require that one updates the state of a numerical model by information from 
sparse and noisy observations, see, e.g., 
\cite{vanLeeuwen2009,Bocquet2010,Fournier2010,PJvL15}
for recent reviews in geophysics. 
Particle filters (PF) are sequential Monte Carlo algorithms
that solve such a ``data assimilation'' problem
as follows \citep{Doucet2001,Tutorial}.
One draws samples from a proposal density,
and then attaches weights to the samples
such that the weighted samples form an empirical 
estimate of a posterior density.
The samples are often called ``particles'' in the PF literature, 
or ``ensemble members'' in the geophysics literature.
The empirical estimate approximates the posterior density in the sense 
that weighted averages over the samples
converge to expected values with respect to the posterior density.
The PFs in the literature differ in their choice of proposal density,
and, therefore in their weights,
see, e.g., 
\cite{OptimalImportanceFunction,GordonSIR,liuchen1995,Doucet1998,Doucet,
Chorin2009,Weare2009,Chorin2010,Morzfeld2011,Weare2012} 
for various particle filters developed over the years.

The variance of the weights determines the efficiency of a PF.
If this variance is large, then a PF ``collapses'', 
i.e., the weights of all but one sample are near zero.
The collapse of particle filters has been described in \cite{Bickel,Bickel2,Bickel3,Snyder,CM13,SBM2015}.
It is argued that an effective dimension defines the ensemble size required to prevent collapse.
While there are various definitions of effective dimensions,
we explain in section~\ref{sec:CollapseReview} that the huge state dimension of global weather models (hundreds of millions of state variables)
and the vast number of observations (millions) cause any effective dimension to be large, 
which in turn implies that particle filters would require huge ensemble sizes in meteorological problems.

It is interesting however that the ensemble Kalman filter (EnKF, \cite{Tippet2003,EvensenBook}),
used routinely with ensembles of 50 to a few hundred members in meteorological problems,
can in fact be interpreted as a particle filter. 
Specifically, as shown by \cite{Papadakis2010}, 
one can view the EnKF ensemble members as draws from a proposal density
and attach weights to each ensemble member,
interpreting the EnKF as a particle filter.
It follows from the optimality result of \cite{SBM2015},
that the particle filter interpretation of the EnKF collapses in the same way as other particle filters.
The goal of this paper is to clarify this situation: 
\emph{Why does the EnKF do so well in meteorological problems
in which its particle filter interpretation collapses?}
Resolving this apparent contradiction 
reveals the significance of a certain type of particle filter collapse,
and suggests improvements of particle filters for meteorological applications.

\section{Preliminaries}
\label{sec:Prelims}
We consider linear Gaussian data assimilation problems,
because we can make our points in this simplified setting.
Moreover, linear problems are standard in the context of the collapse of particle filters 
and have been used by \cite{Bickel,Bickel2,Bickel3,Snyder,CM13,SBM2015}.
Specifically, we consider
\begin{linenomath}
\begin{align}
\label{eq:HMM1}
	\mathbf{x}_k & =\mathbf{M}\mathbf{x}_{k-1}+\mathbf{\eta}_k,\\
\label{eq:HMM2}
	\mathbf{y}_k & = \mathbf{H}\mathbf{x}_k + \mathbf{\varepsilon}_k,
\end{align}
\end{linenomath}
where the state, $\mathbf{x}$,
is a real $n_x$-dimensional vector,
the linear model, $\mathbf{M}$, is an $n_x\times n_x$ matrix,
the observation, $\mathbf{y}$, is a real $n_y$-dimensional vector,
$\mathbf{H}$ is a $n_y\times n_x$ matrix,
and where $k=1,2,3,\dots$ is discrete time.
We further assume that the initial condition $\mathbf{x}_0$ is Gaussian,
and that no observations are available at time $k=0$.
Throughout this paper we use subscripts for time indexing and 
superscripts for indexing samples of random variables.
For example, $\mathbf{x}_k^j$ is the $j$th sample of the state $\mathbf{x}$ at time $k$.
We write $\mathbf{y}_{1:k}$ for 
a set of vectors $\{\mathbf{y}_1,\dots,\mathbf{y}_k\}$.
The random variables $\mathbf{\eta}_k$ and $\varepsilon_k$ 
in equations (\ref{eq:HMM1}) and (\ref{eq:HMM2})
represent model errors and errors in the observations respectively,
and we assume that $\mathbf{\eta}_k$ and $\mathbf{\varepsilon}_k$
are independent identically distributed Gaussians with mean zero,
which we write as
$ \mathbf{\eta}_k\sim\mathcal{N}(\mathbf{0},\mathbf{Q})$,
$ \varepsilon_k\sim\mathcal{N}(\mathbf{0},\mathbf{R})$.
The covariances $\mathbf{Q}$ and $\mathbf{R}$ are $n_x\times n_x$ respectively
$n_y\times n_y$ symmetric positive definite (SPD) matrices.
We exclude data assimilation problems
with a deterministic model for which $\mathbf{Q}=\mathbf{0}$,
however many of our results hold for data assimilation problems with
deterministic models as well by simply setting $\mathbf{Q}=\mathbf{0}$ in the equations. 
We explain which results generalize in section~\ref{sec:Conclusions}.

\subsection{Filtering and smoothing densities}
The goal of data assimilation is to update the state 
of the model~(\ref{eq:HMM1})
in view of the noisy observations in~(\ref{eq:HMM2}).
Model and observations jointly define  
the conditional random variable $\mathbf{x}_{0:k}\vert \mathbf{y}_{1:k}$, 
which describes the state trajectory up to time $k$,
given the observations up to time $k$.
This random variable is specified by the ``smoothing density''
\begin{linenomath}\begin{equation}
  \label{eq:Joint}
  p(\mathbf{x}_{0:k}\vert \mathbf{y}_{1:k}).
\end{equation}\end{linenomath}
The ``filtering density'' is a marginal of the smoothing density,
\begin{linenomath}\begin{equation}
  \label{eq:Marginal}
  p(\mathbf{x}_k\vert \mathbf{y}_{1:k})=\int  \dots \int p(\mathbf{x}_{0:k}\vert \mathbf{y}_{1:k})\, d\mathbf{x}_0 \dots d\mathbf{x}_{k-1}.
\end{equation}\end{linenomath}
and describes the state at observation time, conditioned on all observations
up to now, however it does not describe how the system got to that state.
In practice, the filtering density may be more important.
Forecast models, for example, require initial conditions sampled from the filtering density;
samples from the smoothing density, i.e., descriptions of past weather, are not required.

\subsection{The localized ensemble Kalman filter}
\label{sec:EnKF}
The EnKF uses a Monte Carlo step 
within the Kalman filter formalism,
and recursively approximates the filtering density by an ensemble
\citep{Kalman1960,Kalman1961,Tippet2003,EvensenBook}. 
Specifically, let $\mathbf{x}^j_{k-1}$, $j=1,\dots,N_e$,
be samples of the filtering density $p(\mathbf{x}_{k-1}\vert \mathbf{y}_{1:k-1})$ at time $k-1$. 
The collection of $N_e$ samples is called the ensemble,
and each sample is an ensemble member; $N_e$ is the ensemble size. 
The ensemble is evolved to time $k$ by using the model~(\ref{eq:HMM1})
to generate the forecast ensemble
\begin{linenomath}\begin{equation*}
	\hat{\mathbf{x}}^j_k=\mathbf{M}\mathbf{x}^j_{k-1}+\mathbf{\eta}^j_k,
\end{equation*}\end{linenomath}
where $\mathbf{\eta}^j_k$ is a sample of $\eta_k$. 
Let 
\begin{linenomath}\begin{equation}
	\mathbf{P}^f = \frac{1}{N_e-1}\sum_{j=1}^{N_e} (\mathbf{x}_k^j-\bar{\mathbf{x}}_k) (\mathbf{x}_k^j-\bar{\mathbf{x}}_k)^T,
\end{equation}\end{linenomath}
be the forecast covariance,
where $\bar{\mathbf{x}}_k =\left(\sum \mathbf{x}_k^j\right)/N_e$ is the forecast mean.
The forecast covariance is usually localized and inflated, see, e.g., \cite{HoutMitch2001, Hamilletal2001}.
During localization one decreases long-range correlations in the error covariances,
typically by Schur (entrywise) multiplication of the forecast covariance with a specified $n_x\times n_x$ localization matrix $\rho$:
\begin{linenomath}\begin{equation*}
	\mathbf{P}^f \leftarrow \rho \circ \mathbf{P}^f.
\end{equation*}\end{linenomath}
In atmospheric applications, localization matrices as described by \cite{GC99} are commonly used.
During inflation one increases the localized forecast covariance,
e.g., by multiplying the forecast covariances by a number slightly larger than one,
e.g., $\mathbf{P}^f \leftarrow (1+\alpha) \mathbf{P}^f$, $\alpha >0$.
Combined, localization and inflation are known to be key factors in making
the EnKF a practical tool for meteorological data assimilation,
see, e.g., \cite{Houtekamer2005, Wang2007, Reader2012}.
The localized and inflated forecast covariance defines the Kalman gain
\begin{linenomath}\begin{equation}
\label{eq:KF_K}
\mathbf{K} = \mathbf{P}^f\mathbf{H}^T(\mathbf{H}\mathbf{P}^f\mathbf{H}^T+\mathbf{R})^{-1},
\end{equation}\end{linenomath}
which is used to compute the analysis ensemble
\begin{linenomath}\begin{equation}
  \label{eq:EnKF_Sample}
  \mathbf{x}_k^j = \hat{\mathbf{x}}_k^j+\mathbf{K}\left(\hat{\mathbf{y}}_k^j-\mathbf{H}\hat{\mathbf{x}}_k^j\right),
\end{equation}\end{linenomath}
where $\hat{\mathbf{y}}_k^j =\mathbf{y}_{k}+\varepsilon_k^j$,
is a perturbed observation and where $\varepsilon_k^j$ is a sample of $\varepsilon_k$. 
The sample mean and sample covariance of the analysis ensemble
approximate the mean and covariance of the filtering density $p(\mathbf{x}_k\vert \mathbf{y}_{1:k})$.
This is the ``perturbed observations'' implementation of the localized EnKF.
Other implementations of the EnKF are also available,
see, e.g., \cite{Anderson2001, Bishop2001, Tippet2003}.

\subsection{Particle filters}
Particle filters are sequential importance sampling methods that
construct an empirical estimate of the smoothing density~(\ref{eq:Joint}),
see, e.g., \cite{Doucet2001}.
The empirical estimate consists of a collection of weighted samples,
such that expected values of functions of $\mathbf{x}_{0:k}\vert \mathbf{y}_{1:k}$ can be approximated
by weighted averaging over the samples, see, e.g., \cite{Kalos1986,ChorinHald2013}.
In particular, particle filters make use of a recursion of the smoothing density, 
\begin{linenomath}\begin{equation}
  \label{eq:JointPDFRecursive}
  p(\mathbf{x}_{0:k}\vert \mathbf{y}_{1:k})= p(\mathbf{x}_{0:k-1}\vert \mathbf{y}_{1:k-1})\,\frac{p(\mathbf{x}_{k}\vert \mathbf{x}_{k-1})p(\mathbf{y}_k\vert \mathbf{x}_k)}{p(\mathbf{y}_k\vert \mathbf{y}_{1:k-1})},
\end{equation}\end{linenomath}
and draw samples from a proposal density which is a similar product of functions,
\begin{linenomath}\begin{equation}
  \label{eq:JointImportanceRecursive}
  \pi(\mathbf{x}_{0:k};\mathbf{y}_{1:k})
   \propto \pi_{0}(\mathbf{x}_{0})\prod_{l=1}^{k}\pi_k(\mathbf{x}_l\vert \mathbf{y}_{1:l},\mathbf{x}_{0:l-1}).
\end{equation}\end{linenomath}
Thus,
only the most recent product needs to be updated
when a new observation is collected.  
At each step $k$, the update
$\pi_k(\mathbf{x}_k;\mathbf{y}_{1:k},\mathbf{x}_{0:k-1})$ is used to 
propose samples $\mathbf{x}_k^j$, 
which are appended to the samples of $\mathbf{x}_{0:k-1}^j$
to extend the trajectory to time $k$.
The weights satisfy the recursion
\begin{linenomath}\begin{equation}
\label{eq:weights}
  w_k\propto \frac{q(\mathbf{x}_{0:k}\vert \mathbf{y}_{1:k})}{\pi(\mathbf{x}_{0:k}\vert \mathbf{y}_{1:k})}
  \propto w_{k-1}\frac{p(\mathbf{x}_k\vert \mathbf{x}_{k-1})p(\mathbf{y}_k\vert \mathbf{x}_k)}{\pi_{k}(\mathbf{x}_k\vert \mathbf{y}_{1:k},\mathbf{x}_{0:k-1})},
\end{equation}\end{linenomath}
and account for the fact that one draws samples sequentially from the importance function $\pi(\mathbf{x}_{0:k}\vert \mathbf{y}_{1:k})$,
rather than directly from the smoothing density.
After each step, the weights are normalized such that their sum is one,
and the particles are ``resampled'' \citep{Tutorial},
i.e., particles with a small weight are deleted,
and particles with a large weight are repeated.
Various choices for the updates $\pi_k(\mathbf{x}_k\vert \mathbf{y}_{1:k},\mathbf{x}_{0:k-1})$
lead to the various particle filters developed over the past years,
each with different weights, e.g., 
\cite{GordonSIR,Doucet1998,Chorin2009,Weare2009,Chorin2010,Morzfeld2011,Weare2012}.
The optimal particle filter of \cite{OptimalImportanceFunction,liuchen1995,Doucet}
minimizes the variance of the weights, as was shown by \cite{SBM2015}.

\subsection{The collapse of particle filters and effective dimensions}
\label{sec:CollapseReview}
A particle filter collapses if only one particle carries a significant weight.
The collapse occurs when the variance of the weights is large,
and avoiding the collapse can require huge ensemble sizes.
Indeed, it is often said that particle filters require ensemble sizes that are exponential in the ``dimension''.
However, the situation is more delicate.
It was argued by \cite{CM13} that here are several mechanisms that can cause particle filters to collapse,
not merely a large state dimension.
These mechanisms are described by the theory of the collapse of particle filters and various effective dimensions.
We briefly review this theory here because it is needed to identify the cause for 
the collapse of particle filters in meteorological problems.

The collapse of the ``standard particle filter'' with proposal density
$\pi_k = p(\mathbf{x}_k\vert \mathbf{x}_{k-1})$ is described by \cite{Bickel,Bickel2,Bickel3}.
Specifically, in \cite{Bickel}, it is shown that the collapse of this particle filter depends on the variance
\begin{linenomath}\begin{equation*}
\label{eq:TauSquare}
	\tau^2 = var\left(\log w\right) = \sum_{j=1}^{n_y}\frac{1}{2} \lambda_j^2+\lambda_j\,(\mathbf{y}_k)_j^2,
\end{equation*}\end{linenomath}
where $(\mathbf{y}_k)_j$ is the $j$th element of the observation $\mathbf{y}_k$, 
and where $\lambda_j$ are the eigenvalues of the covariance $cov(\mathbf{R}^{-1/2}\mathbf{H}\mathbf{x}_{k})$.
Precise statements about the collapse of the standard particle filter can be made in the limit
of infinitely many dimensions, $n_y\to\infty$, and provided that $\sum\lambda_j\to\infty$ in this limit. 
In this case, collapse occurs if  $(\log N_e)/\tau^2\to 0$,
i.e., the ensemble size must be exponential in $\tau^2$ to avoid the collapse.
Similar arguments hold for the optimal particle filter,
with a slightly different definition of $\tau^2$ \citep{Snyder, SBM2015}.
The exponential scaling of the ensemble size $N_e$ with $\tau^2$, 
suggests that the ensembles size required to avoid collapse may be too large
to be feasible in numerical weather prediction,
where $n_y$ is so large -- typically millions -- that even moderate eigenvalues $\lambda_j$ lead to large $\tau^2$.

Effective dimensions are quantities related to $\tau^2$,
which describe how difficult it is to solve
a given data assimilation problem by a particle filter.
There are several definitions of effective dimensions in the literature,
some of which are reviewed by \cite{Stuart16}.
For example, in \cite{Bickel2} an effective dimension is defined by 
\begin{linenomath}\begin{equation}
\label{eq:TauSPF}
	\tau_\text{spf} = \text{trace}(cov(\mathbf{R}^{-\frac{1}{2}}\mathbf{H}\mathbf{x}_{k})) = \sum_{j=1}^{n_y} \lambda_j.
\end{equation}\end{linenomath}
Note that this effective dimension is  a property of the data assimilation problem \emph{and} the standard particle filter algorithm.
In \cite{CM13}, effective dimensions for the standard and optimal particle filter are defined by
Frobenius norms, rather than traces of covariances. 
The Frobenius norm of a square matrix is the square root of the sum of the squares of the eigenvalues. 
For example, an effective dimension of the standard particle filter is defined by
\begin{linenomath}\begin{equation*}
	\hat{\tau}_\text{spf}^2 = \vert\vert cov(\mathbf{R}^{-\frac{1}{2}}\mathbf{H}\mathbf{x}_{k})\vert\vert_F^2 = 
	\sum_{j=1}^{n_y} \lambda_j^2.
\end{equation*}\end{linenomath}
However, all effective dimensions above depend on the covariance of $\mathbf{x}_{k}$,
which can be computed from a steady state covariance matrix $\mathbf{P}$, see, e.g., \cite{CM13}.
A steady state requires that the system be linear, 
however the trace of the prior and posterior covariance matrices  
often oscillate about a steady-state value in nonlinear problems as well. 
Indeed, the Frobenius norm, or the trace, of a steady state covariance matrix $\mathbf{P}$ is a natural
candidate for an effective dimension of a data assimilation problem,
independently of the particle filter one applies to solve~it:
\begin{linenomath}\begin{equation}
	\label{eq:TauSys}
	\tau_\text{sys}^2 = \vert\vert\mathbf{P}\vert\vert_F^2 = 
	\sum_{j=1}^{n_x} \alpha_j^2. 
\end{equation}\end{linenomath}
Here, $\alpha_j$ are the eigenvalues of  $\mathbf{P}$.
\cite{CM13} argue that one should only attempt to solve problems
with moderate effective dimension $\tau_\text{sys}$.

It was pointed out in \cite{CLMMT15}, that none of the effective dimensions are in fact ``dimensions''.
The effective dimensions above rather describe ``effective variance'', 
and one can define more effective dimensions, for example, by 
\begin{linenomath}\begin{equation}
	 \label{eq:ellsys}
	 \ell_\text{sys} = \min\ell \in \mathbb{Z} : \sum_{j=1}^{\ell} \alpha_j^2 \ge (1-\kappa) \sum_{j=1}^{n_x} \alpha_j^2,
\end{equation}\end{linenomath}
where $\kappa$ is a predetermined small parameter \cite{CLMMT15}. 
Similar ideas were used in the context of ``partial noise'' by \cite{Morzfeld2012}.
However, it is not yet clear which effective dimension will ultimately be most useful.
Nonetheless, all effective dimensions have in common that they describe 
that ``if the effective dimension is large, an even larger ensemble size is needed to avoid collapse''.

\subsection{The collapse of particle filters in meteorology}
\label{sec:CollapseMeteorology}
We are interested in the collapse of particle filters as it may occur in \emph{meteorological problems}.
Meteorological problems are characterized by
\vspace{-2mm}
\begin{enumerate}[(a)]
\item
large state dimension $n_x$ (hundreds of millions), and large observation dimension $n_y$ (millions);
\item 
large spatial domains, in sense that there are many correlation lengths of a prior distribution contained within the domain;
\vspace{-2mm}
\item 
the matrix $\mathbf{H}$ is ``local'', in the sense that each element of the vector $\mathbf{Hx}$ 
depends only in ``a few'' elements of $\mathbf{x}$.
\end{enumerate}
\vspace{-2mm}
These criteria imply that a standard particle filter avoids the collapse only if 
the ensemble size $N_e$ is exponentially related to $\tau^2$ in (\ref{eq:TauSquare}) 
\citep{Bickel}.
Similar conditions hold for the optimal particle filter as was shown by \cite{Snyder, SBM2015,CM13}.
Moreover, under the above assumptions, $\tau^2$ and in fact \emph{all} effective dimensions 
recited in section~\ref{sec:CollapseReview} are large.
The reasons are as follows.
Recall that the Frobenius norm of a matrix is equal to the square root of the sum of
the squares of all of its elements, $\vert\vert \mathbf{P}\vert\vert^2= \sum_{i,j}P_{i,j}^2$.
Thus, if $n_x$ is huge, 
the effective dimension $\tau_\text{sys}$ is large
unless each entry of $\mathbf{P}$ is tiny,
even if the eigenvalues $\alpha_j$  are moderate.
A large $\tau_\text{sys}$ implies that the effective dimensions 
$\tau_\text{spf}$ and $\hat{\tau}_\text{spf}$ are also large, as shown by \cite{CM13}.
Similar arguments apply to effective dimensions defined for the optimal particle filter
\citep{CM13,Snyder,SBM2015}.
Property~(b) implies that not all variables are correlated,
so that the overall system can be thought of as a collection of loosely coupled lower-dimensional ``sub-systems''.
Since the correlation lengths are small compared to the domain,
one can expect that the overall system contains a large number of such sub-systems.
Thus, the effective dimension $\ell_\text{sys}$ can be expected to be large,
and increasing with the number of loosely coupled subsystems.
Moreover, property (b) implies that the covariance $\mathbf{P}$ is sparse.
The effective dimension $\tau_\text{sys}$ is thus smaller than if all variables were coupled,
which corresponds to a dense matrix $\mathbf{P}$,
however $\tau_\text{sys}$ can still be expected to be large in meteorology.
As an extreme example of sparsity, suppose that $\mathbf{P}$ is diagonal,
and that each diagonal element were $O(1)$;
the effective dimension $\tau_\text{sys}$ would be $O(\sqrt{n_x})$,
which, for a system with hundreds of millions of state variables, is about $10^4$.
This is a perhaps crude lower bound for the effective dimension $\tau_\text{sys}$
in meteorology, since all correlations are neglected in this estimate.
Nonetheless, it shows that $\tau_\text{sys}$,
and in fact all effective dimensions, can be expected to be large in meteorology.
Since it was shown by \cite{CM13,SBM2015} that even the optimal particle filter
requires large ensemble sizes if effective dimensions are large,
one can expect that \emph{all} particle filters collapse in meteorological problems unless the ensemble size is huge.
Note that the ``type'' of filter collapse we describe here
occurs in the canonical linear example of \cite{Bickel,Bickel2,Bickel3,Snyder,CM13,SBM2015},
in which the required ensemble size scales exponentially with the state dimension $n_x$
(see section \ref{sec:MSE} below).

\section{The collapse of the particle filter interpretation of the EnKF}
\label{sec:Collapse}
The EnKF performs well in high-dimensional meteorological applications
in which particle filters are expected to collapse.
On the other hand, it was shown by \cite{Papadakis2010} that
the EnKF can be interpreted as a particle filter (PF-EnKF),
and, thus, can be expected to collapse under the same conditions as other particle filters.
We study this apparent contradiction and in particular
the mechanisms that cause the collapse of the PF-EnFK in meteorological problems.
This approach allows us to combine what we know about the performance of EnKF in \emph{real} meteorological problems
with what we know about the collapse of particle filters.
In this context, we first consider the PF-EnKF as a particle filter for the smoothing density,
and then make connections to the filtering density using ideas similar to marginal particle filters introduced by \cite{Klaas}.

\subsection{Collapse of the PF-EnKF for the smoothing density}
\label{sec:Collapse1}
It was pointed out by \cite{Papadakis2010} that the EnKF ensemble members
can be viewed as draws from the proposal density 
\begin{linenomath}\begin{equation}
\label{eq:EnKFImportanceFunction}
  \pi_{k,\text{EnKF}}^{j}(\mathbf{x}_k\vert \mathbf{x}^j_{k-1},\mathbf{y}_k) = \mathcal{N}\left(\mu^j_k,\Sigma_k\right),
\end{equation}\end{linenomath}
where
\begin{linenomath}\begin{equation*}
  \mu_k^j = (\mathbf{I}-\mathbf{KH})\mathbf{M}\mathbf{x}_{k-1}^j+\mathbf{K}\mathbf{y}_k,\quad \Sigma_k = (\mathbf{I}-\mathbf{K}\mathbf{H})\mathbf{Q}(\mathbf{I}-\mathbf{K}\mathbf{H})^T+\mathbf{K}\mathbf{R}\mathbf{K}^T.
\end{equation*}\end{linenomath}
This proposal density is of the form~(\ref{eq:JointImportanceRecursive})
so that the EnKF can be interpreted as a particle filter.
A related construction using Kernel density estimation at each step is discussed in \cite{Khalil2015}.
The weights that can be attached to the EnKF ensemble are the ratio of the smoothing and proposal densities
\begin{linenomath}\begin{equation}
\label{eq:EnKFWeights}
  w^j_{k,\text{EnKF}} \propto w^{j}_{k-1,\text{ EnKF}}\;\frac{p(\mathbf{x}_k^j\vert \mathbf{x}_{k-1}^j)p(\mathbf{y}_k\vert \mathbf{x}_k^j)}{\pi^j_{k,\text{EnKF}}(\mathbf{x}_k\vert \mathbf{x}^j_{k-1},\mathbf{y}_k)}.
\end{equation}\end{linenomath}
The theory for the collapse of particle filters reviewed in section \ref{sec:CollapseReview} applies to these weights,
and their variance governs the collapse of the PF-EnKF.
This variance is at least as large as the corresponding variance of the optimal PF \citep{SBM2015}.
Thus, the PF-EnKF collapses whenever the optimal PF collapses.
Using the EnKF to define a proposal density therefore can not prevent the collapse 
for small ensemble sizes in meteorological problems.
As a specific example, we will demonstrate in section~\ref{sec:MSE}
that the ensemble size required to avoid the collapse of the PF-EnKF scales exponentially with the dimension $n_x$
in the canonical system studied in \cite{Bickel,Bickel2,Bickel3,Snyder,CM13,SBM2015}.

\subsection{Collapse of the PF-EnKF for the filtering density}
\label{sec:Collapse2}
The particle filter interpretation of the EnKF defines a sequential importance sampler for the smoothing density.
However, the EnKF ensemble can also be viewed as draws from a proposal density for the filtering density.
Indeed, the idea of designing particle filters directly for the filtering density,
rather than the smoothing density,
was explored by \cite{Klaas} via ``marginal particle filters''.
We connect these ideas to the EnKF and study the conditions for the collapse of
the particle filter interpretation of the EnKF for the filtering density.
The reason is that it may be possible that the PF-EnKF collapses as a sampler for the smoothing density,
but performs well as a sampler for the filtering density,
and that sampling the filtering density is sufficient in meteorology.
Below we argue that this is not the case,
i.e., the collapse of the PF-EnKF occurs, under the same conditions,
for filtering and smoothing densities, unless the ensemble size is huge.

As before, we interpret the EnKF ensemble members as draws from a suitable proposal density $\pi_\text{EnKF}(\mathbf{x}_k)$.
We wish to compute the weights, 
\begin{linenomath}\begin{equation}
\label{eq:WeightsEnKF}
	w \propto \frac{p(\mathbf{x}_k\vert \mathbf{y}_{1:k})}{\pi_\text{EnKF}(\mathbf{x}_k)}.
\end{equation}\end{linenomath}
such that the weighted ensemble approximates the filtering density,
not the smoothing density as before.
Note that the EnKF proposal density for the smoothing density in~(\ref{eq:EnKFImportanceFunction})
is conditioned on the previous analysis state
but the EnKF proposal density for the filtering density in~(\ref{eq:WeightsEnKF}) is not.
Thus, the EnKF proposal density for the filtering density is not the same
as the EnKF proposal density for the smoothing density.
The weights~(\ref{eq:WeightsEnKF}), 
are difficult to compute in general because the filtering density and
the EnKF proposal density are not known up to a multiplicative constant 
(see the Appendix and literature about marginal particle filters, e.g., \cite{Klaas}, for more detail).
However, for the linear Gaussian problems we consider,
the filtering density is also Gaussian
and we write $p(\mathbf{x}_k\vert \mathbf{y}_{1:k})\sim\mathcal{N}(\mu_k,\mathbf{P}_k)$.
The proposal density generated by the EnKF ensemble is the Gaussian,
\begin{linenomath}\begin{equation}
\label{eq:EnKFProposalFiltering}
	\pi_\text{EnKF}(\mathbf{x}_k) = \mathcal{N}(\mu_{k,\text{EnKF}},\mathbf{P}_{k,\text{EnKF}}),
\end{equation}\end{linenomath}
where $\mu_{k,\text{EnKF}}$ is the average of the analysis ensemble and where 
\begin{linenomath}\begin{equation}
\label{eq:PEnKF}
	\mathbf{P}_{k,\text{EnKF}}=(\mathbf{I}-\mathbf{KH})\mathbf{MP}_{k-1,\text{EnKF}}\mathbf{M}^T(\mathbf{I}-\mathbf{KH})^T + (\mathbf{I}-\mathbf{KH})\mathbf{Q}(\mathbf{I}-\mathbf{KH})^T +\mathbf{KRK}^T.
\end{equation}\end{linenomath}
Under our assumptions, the weights are
\begin{linenomath}\begin{equation}
\label{eq:EnKFWeightsFilter}
	w\propto \frac{\exp\left(-\frac{1}{2}(\mathbf{x}_k^j-\mu_k)^T\mathbf{P}_k^{-1}(\mathbf{x}_k^j-\mu_k)\right)}
	{\exp\left(-\frac{1}{2}(\mathbf{x}_k^j-\mu_{k,\text{EnKF}})^T\mathbf{P}_{k,\text{EnKF}}^{-1}(\mathbf{x}_k^j-\mu_{k,\text{EnKF}})\right)}.
\end{equation}\end{linenomath}
Asymptotically, for large ensembles, the sample mean converges to the posterior mean and the 
sample covariance converges to the posterior covariance, so that the variance of
the PF-EnKF weights for the filtering density goes to zero asymptotically (see, e.g., \cite{Burgers98}).

We are interested in the collapse of the PF-EnKF with finite ensemble sizes,
and study this situation by making simplifying assumptions.
Specifically, we neglect errors in the ensemble mean,
i.e., we set $\mu_{k,\text{EnKF}}=\mu_k$,
and assume that errors in the ensemble covariance matrix are small, 
i.e., $\mathbf{P}_{k,\text{EnKF}} = (1+\beta)\mathbf{P}_{k}$, where $\beta\geq 0$ is a small number.
While drastic, these simplifications do lead to useful expansions that can accurately describe 
the behavior of the weights of a localized EnKF with small ensembles, 
as will be shown in section~\ref{sec:MSEExample}.
Moreover, our simplifications are ``in favor'' of the EnKF,
i.e., we investigate its performance in a best-case scenario,
that also accounts for the advantageous effects of localization and inflation
since the prior covariance has the right structure (it is a scalar multiple of the posterior covariance), 
and since sampling errors are small.
One could also consider perturbations of the form $(1-\beta)\mathbf{P}_{k}$,
however, in the context of importance sampling, 
perturbations as above are more favorable.
The reason is that importance sampling requires that the support of the proposal density
contains the support of the target density, 
i.e., one wants to use proposal densities with variances larger than those of the target densities.

Under our assumptions,
and if all eigenvalues of the posterior covariance $\mathbf{P}_k$ are positive,
the weights in equation~(\ref{eq:EnKFWeightsFilter}) become
\begin{linenomath}\begin{equation}
\label{eq:WeightsSimple}
	w\propto \exp\left(\frac{1}{2}\mathbf{s}^T\mathbf{s}\right), \quad 
	\mathbf{s} = \sqrt{\frac{\beta}{1+\beta}}\,\left(\mathbf{P}_k^{-1/2}(\mathbf{x}_k^j-\mu_k)\right),
\end{equation}\end{linenomath}
where $\mathbf{x}_k^j\sim\mathcal{N}(\mathbf{\mu}_k,\mathbf{P}_{k,\text{EnKF}})$,
$\mathbf{s} =\mathcal{N}(\mathbf{0},\beta\mathbf{I})$.
The variance of these weights can be computed from the expected values
\begin{align*}
	E\left[w\right] &= \int_{\infty}^\infty \exp(-\mathbf{s}^T\mathbf{s}/2)\frac{\exp(-\mathbf{s}^T\mathbf{s}/(2\beta))}{(2\pi \beta)^{n/2}}  \,\text{d}\mathbf{s} = \left(\frac{1}{1+\beta}\right)^{{n_x}/2},\\
	E\left[w^2\right] &= \int_{\infty}^\infty \exp(-\mathbf{s}^T\mathbf{s})\frac{\exp(-\mathbf{s}^T\mathbf{s}/(2\beta))}{(2\pi \beta)^{n/2}}  \,\text{d}\mathbf{s} = \left(\frac{1}{1+2\beta}\right)^{{n_x}/2}.
\end{align*}
Indeed, one can define an ``effective number of samples'' based on the variance of the weights \citep{Tutorial,Owen}
\begin{linenomath}\begin{equation}
\label{eq:G}
	 N_{e,\text{eff}} = \frac{N_e}{G},\quad G =\frac{var(w)}{E(w)^2}-1=\frac{E(w^2)}{E(w)^2}.
\end{equation}\end{linenomath}
Thus, the collapse occurs when $N_{e,\text{eff}}=1$, i.e., if $G=N_e$,
and the ``quality measure'' $G$ can be used to estimate the ensemble size required to avoid collapse.
Under our assumptions,
\begin{linenomath}\begin{equation*}
	G(\beta) = \left(\frac{1+\beta}{\sqrt{1+2\beta}}\right)^{n_x}
	\approx 1+\frac{\beta^2\,{n_x}}{2}+O(\beta^3).
\end{equation*}\end{linenomath}
The approximation is obtained by expanding $G(\beta)$ around $\beta=0$,
i.e., for small errors $\beta$ in the covariance matrix.
Indeed, if we interpret $\beta\propto 1/\sqrt{N_e}$ as sampling error,
then $G$ can remain constant as dimension $n_x$ increases, 
provided that the ensemble size increases linearly with the dimension.
Thus, the collapse of the PF-EnKF for the filtering density does not occur
for ensemble sizes proportional to the dimension $n_x$.

We can relax the assumption that all eigenvalues of the posterior
covariance $\mathbf{P}_k$ are positive by making use of the effective dimension
$\ell_\text{sys}$ in equation~(\ref{eq:ellsys}).
If some eigenvalues are zero, or if some are close to zero,
then one can substitute the $\ell_\text{sys}\leq n_x$ 
dimensional approximation of this covariance matrix $\mathbf{P}_k$
into the above formulas.
Using the pseudo-inverse of this matrix in (\ref{eq:WeightsSimple}),
one obtains the same expression for $G$ as above with $\ell_\text{sys}$ replacing $n_x$.
We argued above that while $\ell_\text{sys}$ may be smaller than $n_x$ in meteorological problems,
$\ell_\text{sys}$ can be expected to be large in these problems (see property (b) in section~\ref{sec:CollapseMeteorology}).
A large $\ell_\text{sys}$ implies that the collapse of the particle filter interpretation of the EnKF can only be avoided by large ensemble sizes, 
in particular much larger than the 50-100 ensemble members used in practice.

\section{Performance of the EnKF during the collapse of its particle filter interpretation}
\label{sec:MSE}
We have  shown that the collapse of the PF-EnKF for the filtering density can 
be prevented by ensemble sizes that scale linearly with the dimension $n_x$,
or, more generally, with the effective dimension $\ell_\text{sys}$ in~(\ref{eq:ellsys}).
We further argued that the  large dimensions $n_x$ and $n_y$ in meteorological problems
cause \emph{all} particle filters for the smoothing density, 
including the PF-EnKF, to collapse.
On the other hand, localized EnKFs with small ensemble sizes (50 to several hundred) 
are used in meteorological problems with hundreds of millions of state variables
and millions of observations and produce a mean square error (MSE) comparable to 
the normalized trace of the posterior covariance matrix, see, e.g., \cite{Houtekamer2005, Wang2007, Reader2012}.
We now demonstrate that this situation indeed occurs: localized EnKFs with small ensemble sizes
yield a small MSE while their particle filter interpretation collapses.

\subsection{Scaling of MSE with dimension}
The mean square error (MSE) is defined by
\begin{linenomath}\begin{equation}
\label{eq:MSE}
\text{MSE} = \frac{1}{n} \sum_{i=1}^n \left((\bar{\mathbf{x}}_k)_i - (\mathbf{x}^t_k)_i\right)^2,
\quad \bar{\mathbf{x}}_k = \frac{1}{N_e}\sum_{j=1}^{N_e} \mathbf{x}_k^j.
\end{equation}\end{linenomath}
Here $\mathbf{x}^t_k$ is the ``true'' state of the model at time $k$, $(\mathbf{x}^t)_i$, 
$i=1,\dots n$, are its $n$ components, 
and $\bar{\mathbf{x}}_k$ is the sample mean at time $k$, 
computed from the EnKF analysis ensemble members $\mathbf{x}^j_k$,
$j=1,\dots,N_e$.  
We obtain useful scaling equations under the same assumptions as above,
i.e., we consider a Gaussian filtering density $\mathcal{N}(\mu_k,\mathbf{P}_k)$ 
and assume that the EnKF ensemble members are draws from 
$\mathcal{N}(\mu_k,(1+\beta)\mathbf{P}_k)$, where $\beta >0$.
As indicated above, these assumptions simplify the problem,
however describe a best-case scenario for the EnKF (see section \ref{sec:Collapse2}).

Under our assumptions, the statistics of the ensemble mean are
\begin{linenomath}\begin{equation*}
	\bar{\mathbf{x}}_k \sim \mathcal{N}(\mu_k,\frac{1+\beta}{N_e}\,\mathbf{P}_k).
\end{equation*}\end{linenomath}
If the true state is a sample of the filtering density,
one can combine the above expression for the sample mean with the expression 
(\ref{eq:MSE}) for the MSE to obtain the statistics of the MSE:
\begin{linenomath}\begin{equation}
\label{eq:MSEStats}
	\text{MSE} \sim \mathcal{N}\left(1+\frac{1+\beta}{N_e},
	\frac{2\left(1+N_e+\beta\right)^2}{N_e^2\,n_x}\right).
\end{equation}\end{linenomath}
If we assume that the errors in the covariance matrix are due to sampling error,
we can set $\beta = O(1/\sqrt{N_e})$, so that the mean of the MSE is
\begin{linenomath}\begin{equation}
	E\left(\text{MSE} \right) = 1+O(N_e^{-1})+O(N_e^{-3/2}).
\end{equation}\end{linenomath}
Moreover, the variance of the MSE goes to zero for large $n_x$, independently of $N_e$,
because one averages over $n_x$ dimensions.
In summary, our derivation shows that MSE can be small for moderate
ensemble sizes and huge dimensions $n_x$ and $n_y$,
i.e., in situations where the particle filter interpretation of the EnKF collapses (see section \ref{sec:Collapse2}).

\subsection{Local error and MSE vs. global error and collapse of the PF-EnKF}
Note that MSE describes the error one can expect on average per dimension.
Thus, MSE assesses the performance of the EnKF locally, rather than globally.
This assessment of performance is useful in practice.
For example, errors are usually examined locally,
e.g., one might state that the data assimilation is ``working poorly in the tropics, but well in mid-latitudes''.
Similarly, an EnKF user routinely asserts that data assimilation has succeeded 
if average errors at stations around the globe are small.  

However, when the performance of the PF-EnKF is assessed by the variance of its weights,
then one does the opposite.
The variance of the weights is due to differences between the posterior density and
the proposal density generated by the EnKF ensemble. 
If the proposal and posterior densities are effectively high-dimensional
(with large $l_\text{sys}$ and $\tau_\text{sys}$, see sections~\ref{sec:CollapseReview} and~\ref{sec:CollapseMeteorology}),
the differences, or errors, may be small in each dimension,
corresponding to small $\lambda_j$,
however the errors add up to a large overall error during the weight calculation.
This leads to the collapse of the PF-EnKF,
and happens similarly for filtering or smoothing densities.

We have shown that it can indeed occur that localized EnKFs can produce a small MSE 
while at the same time their particle filter interpretation collapses.
The success of the EnKF in meteorology, where MSE is the gold standard,
thus suggests that assessing errors locally is sufficient for these problems.
There is no need for keeping a ``global'' variance of the weights small:
the EnKF can work well even when this variance is huge,
as is illustrated by the collapse of the particle filter interpretation of the EnKF.
Indeed, this collapse may occur regularly and without significant practical consequences.
However, it is important to note that the particle filter interpretation of the EnKF
performs poorly no matter how we evaluate its performance 
-- by MSE or the variance of the weights.
As we will describe below, this problem is caused by 
the fact that the PF-EnKF is not using localization during the weight calculation,
which is equivalent to its ``update'', or analysis.
More specifically, localization \emph{is} used for generating the samples, 
but it is \emph{not} used when calculating the weights 
(see section~\ref{sec:discussion}).

\subsection{Numerical example}
\label{sec:MSEExample}
We illustrate that the localized EnKF with small ensemble sizes $N_e$ can yield a small MSE 
while at the same time its particle filter interpretation collapses 
by the canonical linear example usually studied in the context of filter collapse
\citep{Bickel,Bickel2,Bickel3,Snyder,CM13,SBM2015}.
The example is defined by
\begin{linenomath}\begin{equation}
 \label{eq:Example}
  \mathbf{M} = \mathbf{I},\; \mathbf{H} = \mathbf{I},\; \mathbf{Q}=\mathbf{I},\; \mathbf{R}=\mathbf{I},
\end{equation}\end{linenomath}
and all effective dimensions reviewed in section~\ref{sec:CollapseReview}
increase with dimension $n=n_x=n_y$.
Thus, for large dimension $n$, any particle filter will collapse
unless the ensemble size is huge.
However, the localized EnKF with small ensemble sizes yields a small MSE
even when its particle filter interpretation collapses.
Here, we demonstrate this scenario by computing 
the quality measure $G$ in~(\ref{eq:G}) for various $n$,
which is equivalent to computing the ensemble size required to avoid collapse.
Specifically, if $G$ depends exponentially on dimension $n$,
then the required ensemble size to avoid collapse also depends exponentially on dimension.
Since all effective dimensions increase with dimension $n$,
the scaling of $G$ also illustrates how the 
required ensemble size scales with effective dimensions.

Computing $G$ for one time step, from $k-1$ to $k$,   
requires that one knows the density of the state at $k-1$.
We assume that this density is Gaussian with
mean $\mu_{k-1}=\mathbf{0}$ and covariance $\mathbf{P}_{k-1}=\mathbf{I}$.
We thus have 
\begin{linenomath}\begin{equation*}
p(\mathbf{x}_k\vert \mathbf{x}_{k-1}^j) = \mathcal{N}(\mathbf{x}^j_{k-1},\mathbf{I}),
\quad p(\mathbf{y}_k\vert \mathbf{x}_k^j)= \mathcal{N}(\mathbf{x}^j_{k},\mathbf{I}),
\quad p(\mathbf{x}_{k-1}\vert \mathbf{y}_{1:k-1})= \mathcal{N}(\mathbf{0},\mathbf{I}).
\end{equation*}\end{linenomath}
We consider the following algorithms:
\begin{enumerate}
\vspace{-1mm}
\item an unlocalized\slash uninflated EnKF and its PF interpretation for the filtering density, $N_e=50$,
\vspace{-3mm}
\item a localized\slash inflated EnKF and its PF interpretation for the filtering density, $N_e=50$,
\vspace{-3mm}
\item a localized\slash inflated EnKF and its PF interpretation for the filtering density, $N_e=n$,
\vspace{-3mm}
\item a localized\slash inflated EnKF and its PF interpretation for the smoothing density, $N_e=50$,
\vspace{-3mm}
\item an optimal PF for the smoothing density with $N_e=50$.
\end{enumerate}
\vspace{-1mm}
For each algorithm, we approximate $G$ in~(\ref{eq:G}) from an ensemble of size $N_e=10^6$
and repeat the calculation 100 times to average out errors in our approximations of $G$.
We do this for problems with dimensions between $n=5$ and $n=500$.

\subsubsection{Implementation, localization and inflation of the EnKF}
We localize the forecast covariance by the identity matrix, 
i.e., we set all off-diagonal elements of the EnKF forecast covariance matrix to zero.  
Note that this is a ``perfect'' localization for this diagonal example and 
therefore we are depicting a best-case scenario for the localized EnKF
and its particle filter interpretation. 
Prior inflation is inversely proportional to the ensemble size,
i.e., we multiply the forecast covariance matrix by $(1+1/N_e)$.
However, we tested our results with and without inflation and 
find that, for this example, the impact of inflation is insignificant compared with localization.  
 
The weights of the PF-EnKF for the smoothing density are computed by~(\ref{eq:EnKFWeights}).
Here we assume that $\mathbf{x}_{k-1}\sim p(\mathbf{x}_{k-1}\vert\mathbf{y}_{1:k-1})$,
i.e., the PF-EnKF starts off with samples from the correct density,
so that $w^{k-1} = 1/N_e$.
The weights of the PF-EnKF for the filtering density are computed by~(\ref{eq:WeightsEnKF}).
The proposal density is given by~(\ref{eq:EnKFProposalFiltering}),
and its covariance is localized by the identity matrix as before.
The filtering density is Gaussian, $p(\mathbf{x}_k\vert \mathbf{y}_{1:k})\sim\mathcal{N}(\mu_k,\mathbf{P}_k)$,
and its mean and covariance are computed by the Kalman filter formulas.

\subsubsection{Results}
The left panels of figure~1 show the values of $G$ for the algorithms for each experiment (dots)
along with the log-linear least-squares approximations and one-standard-deviation confidence intervals.
The right panels show the scaling of MSE with dimension,
specifically the mean and 2 standard deviation error bars of the MSE based on the 100 experiments we perform.
We observe the expected exponential scaling of $G$ for the particle filter interpretation of the localized EnKF for the smoothing density. 
The optimal PF is also characterized by an exponential scaling of $G$, however the rate slower (see also \cite{SBM2015}).
The particle filter interpretations of localized or unlocalized EnKFs for the filtering density 
are also characterized by an exponentially increasing $G$ if the ensemble size is fixed at $N_e=50$.
However, the collapse of the particle filter interpretation of the localized EnKF for the filtering density can be avoided 
if the ensemble size is proportional to dimension, as described by our theoretical considerations above.
Specifically, we observe that the PF-EnKF for the filtering density 
produces the same $G$ independently of $n$ when $N_e= n$.

\begin{figure}[tb]
\centering
\includegraphics[width=.85\textwidth]{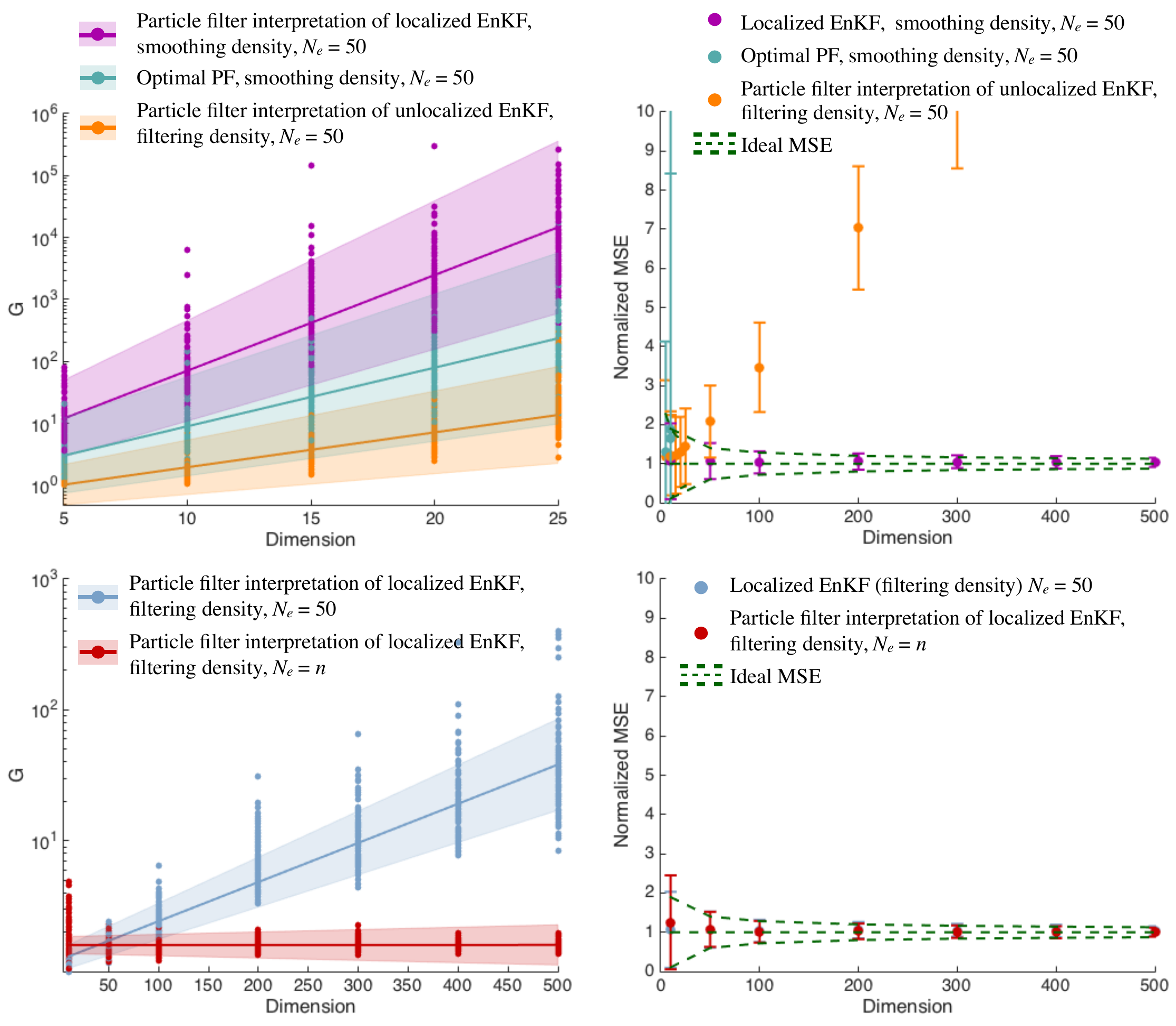}
\caption{
Left: quality measure $G$  as a function of the dimension for 
the various methods. 
Shown are the values of $G$ for each experiment (dots) along with 
log-linear least squares fits (straight lines) and one-standard-deviation confidence intervals (shaded regions).
Right: normalized MSE as a function of the dimension for the various methods.
Shown are the mean (dots) and confidence intervals (error bars) computed from 100 numerical experiments.
Also shown (dashed lines) is the mean and two-standard-deviation confidence intervals of the ideal MSE (dashed green).
}
\label{fig:Collapse}
\end{figure}

In the right panels of figure~1 we observe that the localized EnKFs give a small MSE
while the unlocalized EnKF and the optimal PF cause a large MSE, 
which is consistent with what is reported in practice in large-dimensional problems.
Specifically, the localized EnKF with $N_e=50$ gives a small MSE even 
when its particle filter interpretation collapses.
Also shown is the mean and variance of the MSE 
we obtain without sampling error, i.e., by an optimal algorithm.
The localized EnKFs and the optimal sampling algorithm give almost identical MSE.
In summary, our numerical experiments confirm our scaling of MSE with dimension and 
illustrate that a localized EnKF with a small ensemble size indeed can yield a small MSE, 
while at the same time its particle filter interpretation collapses.

\section{The importance of localization for particle filters}
\label{sec:discussion}
We argued in section~\ref{sec:CollapseMeteorology} that the collapse of particle filters in meteorology is caused by 
extremely large state dimensions $n_x$ and very large numbers of observations $n_y$,
which in turn imply large effective dimensions.
On the other hand, we expect that meteorological problems are characterized by correlation lengths
that are short compared to the domain, and that the observation operator~$\mathbf{H}$ is ``local''.
These assumptions imply that forecast and posterior covariances are sparse,
however the PF-EnKF does not utilize this sparsity during the weight calculation.
If one could make use of the sparsity of the data assimilation problem, 
one may be able to avoid this type of filter collapse of the PF-EnKF, 
even for small ensemble sizes.
The situation is perhaps analogous to numerical linear algebra,
where computations with large dense matrices are expensive,
however large sparse matrices can be dealt with efficiently.

These ideas can be made precise in the diagonal example of section~\ref{sec:MSEExample},
where exploiting sparsity of covariance matrices
means to make use of the fact that the system is diagonal.
Since a diagonal problem uncouples into independent scalar sub-problems, 
one can apply scalar PF-EnKFs to each sub-problem
and also compute the weights for each sub-problem independently.
The ensemble size required for each scalar sub-problem is independent of dimension,
so that the catastrophic scaling of the required ensemble size of the PF-EnKF with dimension disappears.
Indeed, for diagonal problems,
the weight calculation of \emph{every} particle filter can be done in the same way,
i.e., independently for each sub-problem,
so that the scaling of the required ensemble size with dimension disappears for every particle filter
(see also \cite{Rebeschini15}).
These ideas apply equally to particle filters for filtering and smoothing densities,
and are analogous to covariance localization of the EnKF.
Recall that covariance localization of the EnKF reduces or deletes long-range correlations by Schur products (see section~\ref{sec:EnKF}).
The localization thus enforces and exploits sparsity of the forecast covariance,
and in particular restricts the influence of an observation to its neighborhood.
The ``localized'' weight calculation we describe for diagonal problems does precisely that.
Moreover, unlocalized EnKFs with small ensemble sizes fail in high dimensional problems,
and localization makes them effective.
Similarly, particle filters collapse in high-dimensional problems if the ensemble size is small,
however a localization of the weight calcuation makes them effective.

While a localization of the weight calculation for diagonal problems is rather trivial,
the situation is less clear in ``real'' meteorological problems.
In particular, any localization of the weight calculations introduces bias,
as highlighted in the mathematically rigorous paper by \cite{Rebeschini15}.
Localization of forecast covariances of the EnKF also introduces bias.
Moreover, covariance localization of the EnKF and localization of the
weight calculations of particle filters are require the same conditions 
and we argued that these conditions are satisfied in meteorological problems
(see section~\ref{sec:CollapseMeteorology}).
The impressive performance of localized EnKFs in meteorology suggests that the bias 
introduced by the covariance localization is small. 
Thus, if a localization of the weight calculations of particle filters can exploit the same sparsity as covariance localization of the EnKF,
then the introduced bias can also be expected to be small.
Moreover, a localization of the weights of particle filters is only useful if 
local errors such as MSE are of greater importance than 
the global errors measured by the variance of the weights.
We showed in section~{\ref{sec:MSE} 
that the MSE of a localized EnKF with a small ensemble can be small even 
when global errors, measured by the weights of the PF-EnKF, are huge.
The success of localized EnKFs in meteorology thus suggests that it is sufficient to evaluate errors locally.

On the other hand, the success of localization of particle filters also requires
that the relevant effective dimensions of each sub-problem are small enough to prevent the
collapse of particle filters with small ensemble sizes in each sub-problem. 
Our study does not provide estimates of the effective dimension of each sub-problem,
so that, ultimately, we cannot provide proof that localization will indeed make particle filters effective in meteorology.
Our study merely indicates that sparsity of the data assimilation problem can be exploited
by weight localization and that such a strategy can prevent filter collapse in meteorology,
if each sub-problem is manageable.
Indeed, several localization strategies for particle filters have been published recently
\citep{Poterjoy15,Poterjoy16,Penny16,Rebeschini15}.
Earlier work includes the study by \cite{Bengtsson2003},
where localization of mixture filters is considered, 
as well as the paper by \cite{LeiBickel11}  which,
to the best of our knowledge,
is the earliest account of localization for particle filters.
All of these strategies exploit sparsity of the data assimilation problem
by restricting the effects an observation to its neighborhood.
In each paper, simplified numerical examples are provided to demonstrate 
that localized particle filters can solve high dimensional problems with small ensembles.
Our study supports these ideas and connects localization of particle filters to real meteorological problems.
Specifically, we draw attention to the fact that the successful performance of the localized EnKF in numerical
weather prediction in conjunction with the collapse of its particle filter interpretation on the same problem,
indicates that localization schemes for particle filters may prevent particle filter collapse
in meteorology at the cost of a small bias.

\section{Summary}
\label{sec:Conclusions}
We summarize the main results we have collected and derived.
\begin{enumerate}
\item 
We reviewed and extended the literature about the collapse of particle filters
and made explicit connections with the localized EnKF.
Specifically, we showed that the particle filter interpretation of the localized EnKF 
for the filtering density collapses unless the ensemble size increases linearly with dimension.
This means in particular that using the EnKF as the proposal density does not prevent particle filter collapse
in meteorology where the dimension is large.

\item
We have demonstrated that a localized EnKF with a small ensemble size can produce a small MSE even while 
its particle filter interpretation collapses.
The MSE is small because it is a local measure of error,
and the covariance localization of the EnKF can keep local errors small.
The collapse of the particle filter interpretation of the EnKF on the other hand is caused by a global weight calculation,
during which small errors at various locations add up, rather than average out.
We argued that the success of the EnKF indicates that a ``global'' assessment of error may not be required in meteorology
and that assessing errors locally is sufficient.

\item
We argued that the collapse of particle filters can be prevented if 
the data assimilation problem is characterized by sparse posterior and forecast covariances.
This can be done by ``localizing'' the weight calculation,
which is analogous to covariance localization for the EnKF.
Moreover, covariance localization for the EnFK and localization of the weight calculations of particle filters
rely on the same assumptions of short correlation lengths and local assessment of error.
The impressive performance of the localized EnKF
in meteorological problems indicates that these assumptions are indeed valid in meteorology,
which implies that localization of particle filters can be done in meteorology without introducing large bias.
Our collection of results thus corroborates the validity of recent localization schemes
(see, e.g., \cite{PJvL15,Rebeschini15,Poterjoy15,Poterjoy16,Bengtsson2003,Penny16,LeiBickel11}).
\end{enumerate}

The EnKF can struggle in nonlinear problems 
(see, e.g., \cite{Hodyss2011, Hodyss2013})
and in these situations computing weights and using the particle
filter interpretation of the EnKF may be particularly helpful.
While our results do not directly extend to nonlinear problems,
the results we collect for linear problems indicate 
that weight calculations for the EnKF must be localized for 
nonlinear problems as well or else the weights may not be useful in meteorology.

We wish to point out that our results concerning the PF-EnKF for the filtering density in section~\ref{sec:Collapse2},
as well as the results concerning the MSE in section~\ref{sec:MSE} generalize
to data assimilation problems with a ``perfect'' model for which $\mathbf{Q}=\mathbf{0}$.
However, the results concerning the smoothing density cannot be carried over to 
data assimilation problems with deterministic models by setting $\mathbf{Q}=\mathbf{0}$.
The reason is that the smoothing density $p(\mathbf{x}_{1:k}\vert \mathbf{y}_{1:k})$ no longer exists.
In the limit $\mathbf{Q}\to \mathbf{0}$, one can consider either the filtering density
$p(\mathbf{x}_{k}\vert \mathbf{y}_{1:k})$, or the density $p(\mathbf{x}_{0}\vert \mathbf{y}_{1:k})$.
By using the formulas and setting $\mathbf{Q}=\mathbf{0}$, 
one considers methods for the latter density, i.e., 
the density that describes initial conditions conditioned on data collected at later times.
The formulas for the optimal PF with $\mathbf{Q}=\mathbf{0}$ lead to a 
sampling algorithm for $p(\mathbf{x}_{k}\vert \mathbf{y}_{1:k})$, however this algorithm is not optimal.
Other optimal algorithms can be obtained,
as linear examples show.

Finally, we wish to mention that work towards a stability theory
for the EnKF is partially done \citep{Kelly14,Kelly15} and that more is underway \citep{Kelly16a,Kelly16b}.
The goal of these papers is to understand the EnKF and its approximations better
by studying stability of its estimates in \emph{time}.
For example, \cite{Kelly14} interpret the EnKF is as a method 
for approximating the mean of the filtering density,
and the expected distance of an ensemble member to ``the truth'' is analyzed over time. 
It is shown that this distance is bounded as time evolves and can become small,
for large $t$, if the forecast covariances are inflated.
We took a different perspective and address 
scaling of the computational requirements of the particle filter interpretation of the EnKF with dimension.
We considered this scaling during a single time step,
as is customary in the particle filter literature 
\citep{Bickel,Bickel2,Bickel3,Snyder,CM13,SBM2015}.
We thus neglected accumulation of variance and error over time,
and our study complements recent papers about stability of the EnKF, 
as well as earlier results by \cite{DelMoral,Chopin04,Kunsch05},
that address accumulation of variances over time.


\section*{Acknowledgements}
This material is based upon work supported by the 
U.S.~Department of Energy, Office of Science,
Office of Advanced Scientific Computing Research, Applied Mathematics program under contract DE-AC02005CH11231, 
and by the National Science Foundation under grants DMS-1217065 and DMS-1419044.
M.~Morzfeld acknowledges support by an Alfred P.~Sloan Research Fellowship.
D.~Hodyss gratefully acknowledges support from the Office of Naval Research PE-0601153N.
We thank  Prof.~Alexandre J.~Chorin of UC Berkeley and Berkeley National Laboratory 
for interesting discussion and encouragement.

\section*{Appendix: weights for the particle filter interpretation of the EnKF for the filtering density of nonlinear problems}
The weights of the PF-EnKF when sampling the filtering density are given in~(\ref{eq:WeightsEnKF})
and can be re-written as
\begin{linenomath}\begin{equation*}
	w \propto \frac{p(\mathbf{y}_k|\mathbf{x}_k) p(\mathbf{x}_{k}|\mathbf{y}_{1:k-1})}{\pi_\text{EnKF}(\mathbf{x}_k)}.
\end{equation*}\end{linenomath}
Evaluation of these weights requires integration, since
\begin{linenomath}\begin{equation}
\label{eq:Int}
	p(\mathbf{x}_{k}|\mathbf{y}_{1:k-1}) = \int p(\mathbf{x}_{k}|\mathbf{x}_{k-1})p(\mathbf{x}_{k-1}|\mathbf{y}_{1:k-1}) d\mathbf{x}_{k-1},
\end{equation}\end{linenomath}
and evaluating this integral is difficult unless the problem is linear and Gaussian (see above).
Similarly, evaluating the density of the analysis ensemble at time $k$ also requires integration since
\begin{linenomath}\begin{equation}
\label{eq:EnKFUnconditionedPDF}
	\pi_\text{EnKF}(\mathbf{x}_{k}) = \int \pi_{k,\text{EnKF}}(\mathbf{x}_{k}|\mathbf{x}_{k-1},\mathbf{y}_k) p(\mathbf{x}_{k-1}\vert \mathbf{y}_{1:k-1})d\mathbf{x}_{k-1},
\end{equation}\end{linenomath}
where $\pi_{k,\text{EnKF}}(\mathbf{x}_k\vert \mathbf{x}_{k-1},\mathbf{y}_k)$ is not Gaussian. 
However, a basic requirement for weighted sampling (and many other sampling schemes,
e.g., Markov Chain Monte Carlo) 
is that the target and proposal densities be known up to a multiplicative constant.  
This requirement is not met here
so that exact weights for the PF-EnKF when sampling the filtering density are out of reach.

By using ideas of \cite{Klaas} for marginal particle filters,
one can develop approximate weights.
These approximations rely on Monte Carlo for evaluating the integral~(\ref{eq:Int}) 
and become exact as the ensemble size $N_e$ goes to infinity.
We demonstrate here how to use these ideas.
Let 
$\mathbf{x}^j_{k-1}$, $j=1,\dots,N_e$,
be the analysis ensemble at time $k-1$,
which defines an approximation to the filtering density at time $k-1$: 
\begin{linenomath}\begin{equation}
\label{eq:PrevEns}
	p(\mathbf{x}_{k-1}|\mathbf{y}_{1:k-1}) \approx \frac{1}{N_e}\sum_{j=1}^{N_e} \delta(\mathbf{x}_{k-1}-\mathbf{x}^j_{k-1}),
\end{equation}\end{linenomath}
which can be used in the integral~(\ref{eq:Int}) to obtain
\begin{linenomath}\begin{equation*}
	p(\mathbf{x}_{k}|\mathbf{y}_{1:k-1}) \approx \frac{1}{N_e}\sum_{j=1}^{N_e} p(\mathbf{x}_{k}|\mathbf{x}^j_{k-1}).
\end{equation*}\end{linenomath}
The target density can thus be approximated by 
\begin{linenomath}\begin{equation}
\label{eq:MarginalTargetApprox}
	p(\mathbf{x}_{k}|\mathbf{y}_{1:k}) \approx p(\mathbf{y}_k\vert \mathbf{x}_k) \left(\frac{1}{N_e}\sum_{j=1}^{N_e} p(\mathbf{x}_{k}|\mathbf{x}^j_{k-1})\right).
\end{equation}\end{linenomath}
Similarly, we obtain by Monte Carlo approximation of the integral~(\ref{eq:EnKFUnconditionedPDF})
\begin{linenomath}\begin{equation*}
	\pi_\text{EnKF}(\mathbf{x}_k) \approx \frac{1}{N_e}\sum_{j=1}^{N_e}  \pi_\text{EnKF}(\mathbf{x}_k\vert \mathbf{x}^j_{k-1},\mathbf{y}_k).
\end{equation*}\end{linenomath}
With these approximations, the PF-EnKF weights with respect to the marginal posterior are
\begin{linenomath}\begin{equation}
\label{eq:WeightEnKFMarginal}
	w_\text{EnKF}\approx \hat{w}_\text{EnKF} \propto \frac{p(\mathbf{y}_k\vert \mathbf{x}_k)\sum_{j=1}^{N_e} p(\mathbf{x}_k \vert \mathbf{x}^j_{k-1})}{\sum_{j=1}^{N_e}  \pi_\text{EnKF}(\mathbf{x}_k\vert \mathbf{x}^j_{k-1},\mathbf{y}_k)}.
\end{equation}\end{linenomath}

However, the approximate weights are difficult
to compute unless the model's time-stepper is in sync with the observations.
The difficulties can be illustrated by considering a problem with one
additional model step between two consecutive observations.
In this case, the proposal density can be written as 
\begin{linenomath}\begin{equation*}
	\pi_\text{EnKF}(\mathbf{x}_{k}) = \int \int \pi_{k,\text{EnKF}}(\mathbf{x}_{k}\vert \mathbf{x}_{k-1},\mathbf{y}_k) p(\mathbf{x}_{k-1}\vert \mathbf{x}_{k-2})p(\mathbf{x}_{k-2}\vert \mathbf{y}_{1:k-2})d\mathbf{x}_{k-1}d\mathbf{x}_{k-2}.
\end{equation*}\end{linenomath}
We can use the analysis ensemble at time $k-2$ to approximate $p(\mathbf{x}_{k-2}\vert \mathbf{y}_{1:k-2})$ as before,
however we then must perform the integration over $\mathbf{x}_{k-1}$ (assuming there is no observation available at that time):
\begin{linenomath}\begin{equation*}
	\pi_\text{EnKF}(\mathbf{x}_{k}) \approx \frac{1}{N_e}\sum_{j=1}^{N_e} \int \pi_{k,\text{EnKF}}(\mathbf{x}_{k}\vert \mathbf{x}_{k-1},\mathbf{y}_k) p(\mathbf{x}_{k-1}\vert \mathbf{x}^j_{k-2})d\mathbf{x}_{k-1}.
\end{equation*}\end{linenomath}
This integration is difficult in general,
which complicates the computation of the weights for the PF-EnKF
for the filtering density.
In fact, these weights may be impractical,
because the observations are usually not frequently available
(at least not at every step of the model).
These results can be extended to the case of deterministic models
by putting $p(\mathbf{x}_k\vert \mathbf{x}_{k-1})=\delta(\mathbf{x}_k-\mathbf{x}_{k-1})$, where $\delta(\cdot)$ is the delta distribution.

\bibliographystyle{gji}
\bibliography{References}
\end{document}